\documentclass{amsart}
\usepackage{amssymb,amsmath,amsfonts,epsfig,latexsym}

\def\Z{\ensuremath{\mathbb{Z}}} 
\def\P{\ensuremath{\mathbb{P}}}
\def\R{\ensuremath{\mathbb{R}}}

\def\OO{\ensuremath{\mathcal{O}}}

\def\<{\ensuremath{\langle}}
\def\>{\ensuremath{\rangle}}

\begin{document}
    
\title[Toric threefolds with no nontrivial nef line bundles]{Smooth complete toric threefolds with no nontrivial nef line bundles}           
 
\author{Osamu Fujino}
\address{Graduate School of Mathematics, Nagoya University\\ Chikusa-ku\\ Nagoya, 464-8602, Japan.} 
\thanks{The first author is partially supported by The Sumitomo Foundation and by the Grant-in-Aid for Scientific Research \#17684001 from JSPS}    
\author{Sam Payne}
\address{University of Michigan, Department of Mathematics\\ 2074 East Hall, 530 Church St.\\ Ann Arbor, MI 48109, USA.}
\thanks{The second author is supported by a Graduate Research Fellowship from the NSF}

\begin{abstract}
  We describe all of the smooth complete toric threefolds of Picard number 5 with no nontrivial nef line bundles, and show that no such examples exist with Picard number less than 5.
\end{abstract}
\maketitle

\section{Introduction}

\renewcommand{\theenumi}{\alph{enumi}}

Let $X$ be a complete algebraic variety.  A line bundle $L$ on $X$ is
said to be nef if $(L \cdot C)$, the degree of $L$ restricted to $C$, is
nonnegative for every curve $C \subset X$.  In particular, any ample
line bundle on $X$ is nef, as is the trivial line bundle $\OO_{X}$.
Hence complete varieties with no nontrivial nef line bundles are
nonprojective.  Many constructions of smooth or mildly singular
nonprojective varieties use toric geometry; in combinatorial language,
they correspond to nonregular triangulations of totally cyclic vector
configurations, which are also referred to as nonprojective or
noncoherent triangulations.  The standard examples of complete
nonprojective toric varieties \cite{Bonavero}, \cite{Ewald},
\cite[II.2E]{KKMS}, \cite[pp.84--85]{Oda} arise from nonregular
triangulations of boundary complexes of polytopes; such varieties
admit proper birational morphisms to projective toric varietes, and
hence have nontrivial nef line bundles.  If $\Delta$ is the
fan over the faces of a triangulation of the boundary complex of some
polytope $P$, then $X(\Delta)$ is proper and birational over the
projective toric variety $X_{P^{\circ}}$ corresponding to the polar
polytope $P^{\circ}$, and the pullbacks of ample line bundles from
$X_{P^{\circ}}$ are nontrivial nef line bundles on $X(\Delta)$.

If a variety has no nontrivial nef line bundles, then it admits
no nonconstant morphisms to projective varieties.  One example of such a
variety, due to Fulton, is the singular toric threefold corresponding
to the fan over a cube with one ray displaced, which has no nontrivial
line bundles at all \cite[pp.25-26, 72]{Fulton}.  A similar example is
due to Eikelberg \cite[Example 3.5]{Eikelberg}.  The literature also
contains a number of examples of Moishezon spaces of Picard number 1
with no nontrivial nef line bundles.  See \cite{Bonavero},  \cite[5.3.14]{Kollar}, \cite[3.3]{Nakamura},
and \cite{Oguiso}.  Recently,
the first author has given examples of complete, singular toric varieties with
arbitrary Picard number that have no nontrivial nef line bundles \cite{Fujino2}.
Motivated by the observation that all of the standard examples of
smooth, complete, nonprojective toric varieties are proper over
projective varieties, he conjectured that every smooth, complete toric
variety has a nontrivial nef line bundle.  The second author found counterexamples to this conjecture, but the examples were relatively complicated, with Picard number $11$ and higher \cite{Payne}.  In general, it is known that every smooth complete toric variety with Picard number less than 4 is projective \cite{KS}.

This note describes the simplest possible examples of smooth complete toric threefolds with no nontrivial nef line bundles, which have Picard number 5.  These examples consist of two infinite collections of combinatorially equivalent varieties, plus one exceptional example.  The constructions yield similar examples with any Picard number greater than 5, and we show that no such examples exist with Picard number less than 5.  This resolves \cite[Problem 5.2]{Fujino2}, which asked for which Picard numbers there exist smooth complete algebraic varieties with no nontrivial nef line bundles, for the case of toric threefolds.

\vspace{5 pt}

\section{Examples}  We begin by giving examples of smooth toric threefolds with no nontrivial nef line bundles and Picard number 5.

\vspace{5 pt}

\noindent \textbf{Example 1.} Let $\Sigma$ be the fan in $\R^{3}$ whose rays are generated by
\[
\begin{array}{lll}
    v_{1} = (1,0,0), & v_{2} = (0,1,0), & v_{3} = (0,0,1), \\ v_{4} =
    (0,-1,-1), & v_{5} = (-1,0,-1), & v_{6} = (-2,-1,0),
\end{array}
\]
and whose maximal cones are
\[
\begin{array}{llll}
    \< v_{1}, v_{2}, v_{3} \>, & \< v_{1}, v_{2}, v_{4} \>, & 
    \< v_{2}, v_{4}, v_{5} \>, & \< v_{2}, v_{3}, v_{5} \>, \\
    \< v_{3}, v_{5}, v_{6} \>, & \< v_{1}, v_{3}, v_{6} \>, &
    \< v_{1}, v_{4}, v_{6} \>, & \< v_{4}, v_{5}, v_{6} \>.
\end{array}
\]
The nonzero cones of $\Sigma$ are the cones over the faces of the 
following nonconvex polyhedron.

\begin{center}

 \includegraphics{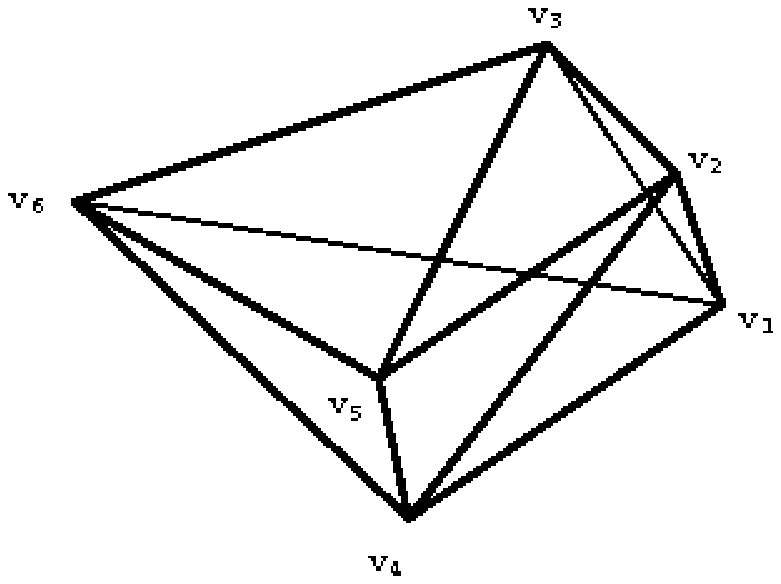}

\end{center}

Let $\Delta$ be the fan obtained from $\Sigma$ by successive star
subdivisions along the rays spanned by $v_7 = (-1,-1,-1)$ and $v_8 = (-2,-1,-1)$.

We claim that $X = X(\Delta)$, the toric threefold corresponding to the
fan $\Delta$ with respect to the lattice $\Z^{3} \subset \R^{3}$, is 
smooth and complete and that $X$ has no nontrivial nef line bundles. 
It is easy to check that $|\Delta| = \R^{3}$, so $X$ is complete, and
that each of the 12 maximal cones of $\Delta$ is spanned by a basis
for $\Z^{3}$, so $X$ is smooth.  It remains to show that $X$ has no
nontrivial nef line bundles.

Suppose $D = \sum d_{i} D_{i}$ is a nef Cartier divisor, where $D_{i}$
is the prime $T$-invariant divisor corresponding to the ray spanned by
$v_{i}$.  Since nef line bundles on toric varieties are globally
generated, $\OO(D)$ has nonzero global sections.  Therefore we may
assume that $D$ is effective, i.e.\ each $d_{i} \geq 0$.  We will show
that $D = 0$.  Note that $\< v_{1}, v_{2}, v_{4} \>$ is a cone in
$\Delta$, and
\[
   v_{5} \, = \, v_{2} + v_{4} - v_{1}.
\]
It then follows from the convexity of the piecewise linear function 
$\Psi_{D}$ associated to $D$ that
\[
   d_{1} + d_{5} \geq d_{2} + d_{4}.
\]   
Similarly,
\[
\begin{array}{lll}
   d_{2} + d_{6} & \geq & 2d_{3} + 2 d_{5}, \\  
   d_{3} + d_{4} & \geq & 2d_{1} + d_{6}.
\end{array}   
\]
Adding the three above inequalities, we have
\begin{eqnarray*}
    \lefteqn{ d_{1} + d_{2} + d_{3} + d_{4} + d_{5} + d_{6} 
    \ \geq } \\
    & & 2d_{1} + d_{2} + 2d_{3} + d_{4} + 2d_{5} + d_{6}.
\end{eqnarray*}  
Since all of the $d_{i}$ are nonnegative, it follows that $d_{1} = d_3 = d_5 = 0.$
Substituting $d_1 = d_5 = 0$ in the first inequality gives $d_2 = d_4 = 0$.
Since $\{v_{1}, \ldots, v_{5}\}$ positively spans
$\R^{3}$, we have $P_D = \{ 0 \}$, where $P_D$ is the polytope defined by
\[
P_D = \{u \in \R^3 : \< u, v_i \> \geq -d_i \mbox{ for all } i \}.
\]
Since $D$ is nef,
\[
-d_i = \min\{ \<u,v_i \> : u \in P_D \},
\]
for all $i$.  Therefore, $D = 0$. 

\vspace{5 pt}

The above argument shows that any complete fan in $\R^{3}$ containing
the cones $\< v_{1}, v_{2}, v_{4} \>,$ $\<v_{2}, v_{3}, v_{5} \>$,
and $\< v_{1}, v_{3}, v_{6} \>$
corresponds to a complete toric variety with no nontrivial nef line bundles.  In
particular, star subdivisions of the remaining cones in $\Delta$
lead to examples of smooth complete toric threefolds with no nontrivial nef line bundles and any Picard 
number greater than or equal to 5.

Example 1 is the toric variety labeled as [8-12] in \cite[p.79]{Oda1}, as can be checked by setting $n = (-1,-1,-1), n' = (1,0,0),$ and $n'' = (0,0,1)$.

\vspace{10 pt}

\noindent \textbf{Example 2.}  Let $a$ be an integer, with $a \neq 0,-1$.  Let $\Delta$ be the fan whose rays are spanned by
\[
\begin{array}{llll}
    v_{1} = (1,0,0), & v_{2} = (0,1,0), & v_{3} = (0,0,1), & v_{4} =
    (0,-1,-a), \\ v_{5} = (0,0,-1), & v_{6} = (-1,1,-1), & 
    v_7 = (-1,0,-1), & v_8 = (-1,-1,0),
\end{array}
\]
and whose maximal cones are
\[
\begin{array}{llll}
    \< v_{1}, v_{2}, v_{3} \>, & \< v_{1}, v_{3}, v_{4} \>, & 
    \< v_{1}, v_{4}, v_{5} \>, & \< v_{1}, v_{5}, v_{6} \>, \\
    \< v_{1}, v_{2}, v_{6} \>, & \< v_{2}, v_{3}, v_{8} \>, &
    \< v_{3}, v_{4}, v_{8} \>, & \< v_{4}, v_{5}, v_{8} \>, \\
    \<v_5, v_6, v_7\>, & \<v_5, v_7, v_8 \>, &
    \<v_6, v_7, v_8\>, & \<v_2, v_6, v_8 \>.
\end{array}
\]
It is straightforward to check that $X(\Delta)$ is smooth and complete.  We claim that $X(\Delta)$ has no nontrivial nef line bundles.  Suppose $D = \sum d_i D_i$ is an effective nef divisor on $X(\Delta)$.  The convexity of $\Psi_D$ gives the following inequalities:
\[
\begin{array}{lll}
   d_{1} + d_{8} & \geq & ad_{3} + d_{4}, \\
   2d_4 + d_6 & \geq  & (2a + 1) d_5 + d_8, \\ 
   d_2 + d_5 & \geq & d_1 + d_6, \\   
   d_{3} + d_{6} & \geq & 2d_{2} + d_{8}, \\
   d_2 + d_4 & \geq & -ad_3,\\
   d_1 + d_8 & \geq & d_4 -ad_5.
\end{array}   
\]
We claim that $d_3 = d_5 = 0$.  When $a >0$, the claim follows from the first four inequalities, since
\begin{eqnarray*}
\lefteqn{2(d_1+d_8) +  2d_4  +  d_6  + 2(d_2+d_5) + d_3 + d_6  \geq}  \\
& & 2(ad_3 + d_4) + (2a+1)d_5 + d_8 + 2(d_1 + d_6) + 2d_2 + d_8.
\end{eqnarray*}
Similarly, the claim follows from the last four inequalities when $a < -1$.  Substituting $0$ for $d_3$ and $d_5$ in the middle two inequalities gives $d_1 = d_2 = d_8 = 0$.  Since $\{v_1, v_2, v_3, v_5, v_8 \}$ positively spans $\R^3$, it follows that $P_D = \{ 0 \}$, and hence $D=0$.  So $X(\Delta)$ has no nontrivial nef line bundles.

\vspace{5 pt}

The varieties in Example 2 are labeled [8-5$'$] in \cite[p.78]{Oda1}, as can be checked by setting $n = (1,0,0), n' = (0,1,0)$, and $n'' = (0,0,1)$.

\vspace{10 pt}

\noindent \textbf{Example 3.}  Let $a$ and $b$ be nonzero integers, with $(a,b) \neq (\pm1, \pm1)$.  Let $\Delta$ be the fan whose rays are spanned by
\[
\begin{array}{llll}
v_1 = (-1,b,0), &  v_2 = (0,-1,0), & v_3 = (1,-1,0), &
v_4 = (-1,0,-1), \\ v_5 = (0,0,-1), & v_6 = (0,1,0), &
v_7 = (0,0,1), & v_8 = (1,0,a),
\end{array}
\]
and whose maximal cones are
\[
\begin{array}{llll}
\<v_1, v_2, v_4\>, & \<v_2, v_3, v_4\>, & \<v_3, v_4, v_5\>, & \<v_4, v_5, v_6 \>,\\
\<v_1, v_4, v_6\>, & \<v_1, v_6, v_7\>, & \<v_1, v_2, v_7\>, & \<v_2, v_3, v_7\>, \\
\<v_3, v_5, v_8\>, & \<v_5, v_6, v_8\>, & \<v_3, v_7, v_8\>, & \<v_6, v_7, v_8\>.
 \end{array}
 \]
 It is straightforward to check that $X(\Delta)$ is smooth and complete.  We claim that $X(\Delta)$ has no nontrivial nef line bundles.  Suppose $D = \sum d_i D_i$ is an effective nef divisor on $X(\Delta)$.  We now consider the case where $a$ and $b$ are both positive, and show that $D = 0$.  The convexity of $\Psi_D$ gives the following inequalities:
\[
 \begin{array}{lll}
 d_3 + d_6 & \geq & ad_5 + d_8, \\
 d_2 + d_8 & \geq & d_3 + ad_7, \\
 d_1 + d_5 & \geq & bd_6 + d_4, \\
 d_4 + d_7 & \geq & d_1 + bd_2. \\
 \end{array}
 \]
Since $(a,b) \neq (1,1)$, either $a>1$ or $b >1$.  We claim that $d_2 = d_5 = d_6 = d_7 = 0$.  If $a > 1$, then adding the above inequalities gives $d_5 = d_7 = 0$, and substituting 0 for $d_5$ and $d_7$ in the last two inequalities gives $d_2 = d_6 = 0$.  If $b>1$, then adding the above inequalities gives $d_2 = d_6 = 0$, and substituting 0 for $d_2$ and $d_6$ in the first two inequalities gives $d_5 = d_7 = 0$.  This proves the claim.
Now, convexity of $\Psi_D$ gives
 \[
 \begin{array}{lll}
 d_2 + d_5 & \geq & d_3 + d_4,
 \end{array}
 \]
 so $d_3 = d_4 = 0$.  Since $\{ v_2, \ldots, v_7 \}$ positively spans $\R^3$, it follows that $P_D = \{ 0 \}$, and hence $D = 0$.
 
The cases where $a$ or $b$ is negative, or both, are proven similarly, using the inequalities above, together with the following inequalities, which also come from the convexity of $\Psi_D$:
 \[
 \begin{array}{lll}
 d_3 + d_6 & \geq & -ad_7 + d_8, \\ 
 d_4 + d_8 & \geq & (1-a)d_5, \\
 d_4 + d_7 & \geq & d_1 - bd_6, \\
 d_1 + d_3 & \geq & (1-b)d_2.
 \end{array}
 \]
  
  The varieties in Example 3 are labeled [8-13$'$] in \cite[p.79]{Oda1}, as can be checked by setting $n = (1,-1,0), n' = (0,1,0),$ and $n'' = (0,0,1)$.

\section{Classification}  Using the classification of minimal smooth toric threefolds with Picard number at most 5, due to Miyake and Oda \cite[Section 9]{Oda1}, and Nagaya \cite{Nagaya}, we show that Examples 1, 2, and 3 give all such varieties that have no nontrivial nef line bundles.  Recall that a smooth toric variety is minimal if it is not the blow up of a smooth toric variety along a smooth torus invariant center.

\vspace{10 pt}

\noindent \textbf{Theorem}  \emph{Every smooth complete toric threefold with no nontrivial nef line bundles and Picard number at most 5 is isomorphic to one described in Examples 1, 2, and 3, above. }

\vspace{10 pt}

It will suffice to show that there are no minimal smooth complete toric threefolds with no nontrivial nef line bundles and Picard number less than 5, and that every minimal smooth complete nonprojective toric threefold with Picard number 5 that is not listed in Examples 1, 2, and 3, admits a nonconstant morphism to a projective variety.

\vspace{10 pt}

\noindent \textbf{Lemma} \emph{There are no minimal smooth complete toric threefolds with no nontrivial nef line bundles and Picard number less than 5.}

\vspace{10 pt}

\noindent \emph{Proof.} By the classification theorem of Miyake, Oda, and Nagaya, every minimal smooth complete toric threefold with Picard number less than 5 is isomorphic to either $\P^3$, a $\P^2$-bundle over $\P^1$, a $\P^1$-bundle over a projective toric surface, or one of the following:

\begin{enumerate} 

\item The projective toric variety $X(\Delta)$, where $\Delta$ is the fan in $\R^3$ whose rays are generated by
\[
\begin{array}{lll}
v_1 = (1,0,0), & v_2 = (0,1,0), &
 v_3 = (0,0,1), \\ v_4 = (0,-1,a), &
  v_5 = (0,0,-1), & v_6 = (0,1,-1), \\
v_7 = (-1,2,-1), & &
\end{array}
\]
for some integer $a$, and whose maximal cones are
\[
\begin{array}{llll}
\<v_1, v_2, v_3\>, & \<v_1, v_3, v_4\>, & \<v_1, v_4, v_5\>, & \<v_1, v_5, v_6\>, \\
\<v_1, v_6, v_7\>, & \<v_1, v_2, v_7\>, & \<v_2, v_3, v_7\>, & \<v_3, v_4, v_7\>, \\
\<v_4, v_5, v_7\>, & \<v_5, v_6, v_7\>. & &
\end{array}
\]

\item The toric variety $X$ associated to the fan in $\R^3$ whose rays are generated by
\[
\begin{array}{lll}
v_1 = (-1,0,0), & v_2 = (0,-1,0), & v_3 = (0,0,-1), \\
v_4 = (1,0,1),  & v_5 = (0,1,1), & 
 v_6 = (1,1,1), \\ v_7 = (1,1,0), & &
\end{array}
\]
and whose maximal cones are
\[
\begin{array}{llll}
\<v_1, v_2, v_3\>, & \<v_2, v_3, v_4\>, & \<v_2, v_4, v_5\>, & \<v_1, v_2, v_5\>, \\
\<v_4, v_5, v_6\>, & \<v_1, v_3, v_7\>, & \<v_3, v_4, v_7\>, & \<v_4, v_6, v_7\>, \\
\<v_5, v_6, v_7\>, & \<v_1, v_5, v_7\>. & &
\end{array}
\]
In this case, $X$ is not projective \cite[Proposition 9.4]{Oda1}, but the anticanonical bundle $-K_X$ is nef.
\end{enumerate}
\vspace{5 pt}
Therefore, up to isomorphism, there is only one nonprojective smooth complete toric threefold with Picard number less than 5, and it has a nontrivial nef line bundle. \hfill $\Box$

\vspace{10 pt}

\noindent \emph{Proof of Theorem.}  It remains to show that every minimal smooth complete nonprojective toric threefold of Picard number 5 that is not isomorphic to one in Examples 1, 2, or 3 has a nontrivial nef line bundle.  By the classification theorem of Miyake and Oda, and Nagaya, and by the second Remark in \cite[p.80]{Oda1}, any such variety is isomorphic to one of the varieties listed below.  In each case, it will suffice to either exhibit a specific nontrivial nef line bundle, or to give a nonconstant morphism to a projective variety.  Note that the varieties $X(\Delta)$ as in Examples 2 and 3, but with $a = 0$ or $b = 0$, are not minimal.  The labels such as [8-$5'$] are the labels given to the varieties in \cite[pp.78-79]{Oda1}.  Different labels with the same numbers, such as [8-5$'$] and [8-5$''$], indicate that the corresponding varieties have combinatorially equivalent fans.

\begin{enumerate}  

\item \mbox{[8-5$'$]} The toric variety $X = X(\Delta)$ as in Example 2, with $a = -1$.  In this case, $-K_X - D_4 - D_8$ is nef.

\item  \mbox{[8-5$''$]} The toric variety $X(\Delta)$ associated to the fan in $\R^3$ whose rays are generated by
\[
\begin{array}{lll}
v_1 = (0,1,0), & v_2 = (0,-1,-1), &
 v_3 = (1,0,0), \\ v_4 = (0,0,1), &
  v_5 = (-1,0,-1), & v_6 = (-1,-2,-2), \\
v_7 = (-1,-1,-1), & v_8 = (-1,-1,0), 
\end{array}
\]
and whose maximal cones are
\[
\begin{array}{llll}
    \< v_{1}, v_{2}, v_{3} \>, & \< v_{1}, v_{3}, v_{4} \>, & 
    \< v_{1}, v_{4}, v_{5} \>, & \< v_{1}, v_{5}, v_{6} \>, \\
    \< v_{1}, v_{2}, v_{6} \>, & \< v_{2}, v_{3}, v_{8} \>, &
    \< v_{3}, v_{4}, v_{8} \>, & \< v_{4}, v_{5}, v_{8} \>, \\
    \<v_5, v_6, v_7\>, & \<v_5, v_7, v_8 \>, &
    \<v_6, v_7, v_8\>, & \<v_2, v_6, v_8 \>.
\end{array}
\]
In this case, $\Delta$ is a refinement of the complete fan whose set of rays is $\{v_1, v_3, v_4, v_7\}$.  Therefore $X(\Delta)$ admits a proper birational morphism to $\P^3$.

\item \mbox{[8-8]} The toric variety $X(\Delta)$ associated to the fan in $\R^3$ whose rays are generated by
\[
\begin{array}{lll}
v_1 = (0,0,1), & v_2 = (1,0,0), &
v_3 = (0,-1,-1), \\ v_4 = (-1,-2,-1), &
v_5 = (0,1,0), & v_6 = (0,0,-1), \\
v_7 = (-1,-2,-2), & v_8 = (-1,-1,-2), 
\end{array}
\] 
and whose maximal cones are
\[
\begin{array}{llll}
\<v_1, v_2, v_3\>, & \<v_1, v_3, v_4\>, & \<v_1, v_4, v_5\>, & \<v_1, v_2, v_5\>, \\
\<v_2, v_5, v_6\>, & \<v_3, v_4, v_7\>, & \<v_2, v_3, v_8\>, & \<v_3, v_7, v_8\>, \\
\<v_4, v_7, v_8\>, & \<v_4, v_5, v_8\>, & \<v_5, v_6, v_8\>, & \<v_2, v_6, v_8\>.
\end{array}
\]
In this case, $\Delta$ is a subdivision of the complete fan whose rays are $v_1, v_2, v_5, v_7$.  Therefore $X(\Delta)$ admits a proper birational morphism to $\P(1,1,2,2)$.

\item \mbox{[8-11]}
The toric variety $X(\Delta)$ associated to the fan in $\R^3$ whose rays are generated by
\[
\begin{array}{lll}
v_1 = (1,0,0), & v_2 = (0,-1,0), &
v_3 = (0,0,1),\\ v_4 = (1,1,a), &
v_5 = (0,0,-1), & v_6 = (0,-1,-1), \\
v_7 = (-1,-2,-1), & v_8 = (0,1,b),
\end{array}
\]
for some integers $a$ and $b$, and whose maximal cones are
\[
\begin{array}{llll}
\<v_1, v_2, v_3 \>, & \<v_1, v_3, v_4 \>, & \<v_1, v_4, v_5 \>, & \<v_1, v_5, v_6\>, \\
\<v_1, v_6, v_7 \>, & \<v_1, v_2, v_7 \>, & \<v_2, v_3, v_7 \>, & \<v_5, v_6, v_7 \>, \\
\<v_3, v_4, v_8 \>, & \<v_4, v_5, v_8 \>, & \<v_5, v_7, v_8 \>, & \<v_3, v_7, v_8\>.
\end{array}
\]
In this case, the projection $(x,y,z) \mapsto (x,y)$ gives a map from $\Delta$ to the complete fan in $\R^2$ whose rays are generated by $(1,0), (0,1), (1,1)$, and $(-1,-2)$.  Therefore $X(\Delta)$ admits a proper surjective morphism to the blow up of $\P(1,1,2)$ at a point.

\item \mbox{[8-13$'$]} The toric variety $X = X(\Delta)$ as in Example 3, with $(a,b) = (\pm 1, \pm 1)$.  In this case, $-K_X$ is nef.

\item \mbox{[8-13$''$]}
The toric variety $X(\Delta)$ associated to the fan in $\R^3$ whose rays are generated by
\[
\begin{array}{lll}
v_1 = (1,1,b), & v_2 = (1,0,0), &
v_3 = (0,-1,0), \\ v_4 = (0,0,-1), &
v_5 = (-1,a,d), & v_6 = (0,1,0), \\
v_7 = (0,0,1), & v_8 = (-1,c,d+1),
\end{array}
\]
for some integers $a, b, c,$ and $d$, and whose maximal cones are
\[
\begin{array}{llll}
\<v_1, v_2, v_4\>, & \<v_2, v_3, v_4\>, & \<v_3, v_4, v_5\>, & \<v_4, v_5, v_6 \>,\\
\<v_1, v_4, v_6\>, & \<v_1, v_6, v_7\>, & \<v_1, v_2, v_7\>, & \<v_2, v_3, v_7\>, \\
\<v_3, v_5, v_8\>, & \<v_5, v_6, v_8\>, & \<v_3, v_7, v_8\>, & \<v_6, v_7, v_8\>.
 \end{array}
 \]
In this case, the projection $(x,y,z) \mapsto x$ gives a map from $\Delta$ to the complete fan in $\R$.  Therefore $X(\Delta)$ admits a proper surjective morphism to $\P^1$.

\item \mbox{[8-14$'$]}
The toric variety $X(\Delta)$ associated to the fan in $\R^3$ whose rays are generated by
\[
\begin{array}{lll}
v_1 = (0,0,-1), & v_2 = (1,0,0), &
v_3 = (0,1,0), \\ v_4 = (-1,-1,a), &
v_5 = (-1,-1,a+1), & v_6 = (1,0,1), \\
v_7 = (0,0,1), & v_8 = (0,1,1),
\end{array}
\]
for some integer $a$, and whose maximal cones are
\[
\begin{array}{llll}
\<v_1, v_2, v_3\>, & \<v_1, v_3, v_4\>, & \<v_1, v_2, v_4\>, & \<v_3, v_4, v_5\>, \\
\<v_4, v_5, v_6\>, & \<v_2, v_4, v_6\>, & \<v_5, v_6, v_7\>, & \<v_2, v_3, v_8\>, \\
\<v_3, v_5, v_8\>, & \<v_5, v_7, v_8\>, & \<v_6, v_7, v_8\>, & \<v_2, v_6, v_8\>.
\end{array}
\]
In this case, the projection $(x,y,z) \mapsto (x,y)$ gives a map from $\Delta$ to the complete fan in $\R^2$ whose rays are generated by $(1,0), (0,1),$ and $(-1,-1)$.  Therefore $X(\Delta)$ admits a proper surjective morphism to $\P^2$.

\item \mbox{[8-14$''$]}
The toric variety $X(\Delta)$ associated to the fan in $\R^3$ whose rays are generated by
\[
\begin{array}{lll}
v_1 = (-1,a,b), & v_2 = (0,1,0), &
v_3 = (0,-1,-1), \\ v_4 = (0,0,1), &
v_5 = (1,0,1), & v_6 = (1,1,0), \\
v_7 = (1,0,0), & v_8 = (1,-1,-1),
\end{array}
\]
for some integers $a$ and $b$, and whose maximal cones are
\[
\begin{array}{llll}
\<v_1, v_2, v_3\>, & \<v_1, v_3, v_4\>, & \<v_1, v_2, v_4\>, & \<v_3, v_4, v_5\>, \\
\<v_4, v_5, v_6\>, & \<v_2, v_4, v_6\>, & \<v_5, v_6, v_7\>, & \<v_2, v_3, v_8\>, \\
\<v_3, v_5, v_8\>, & \<v_5, v_7, v_8\>, & \<v_6, v_7, v_8\>, & \<v_2, v_6, v_8\>.
\end{array}
\]
In this case, the projection $(x,y,z) \mapsto x$ gives a map from $\Delta$ to the complete fan in $\R$.  Therefore $X(\Delta)$ admits a proper surjective morphism to $\P^1$.
\end{enumerate}
\vspace{-13 pt} \hfill $\Box$

\end{document}